\documentclass[10pt,a4paper,reqno]{amsart} 
\usepackage{ae,aecompl}
\usepackage[cm]{aeguill}
\usepackage{pstricks,amsmath,here,latexsym,amssymb}
\usepackage{fullpage} 
\usepackage{epsf,psfrag,psfig,epsfig,color,graphicx}
\usepackage[latin1]{inputenc}
\usepackage[english]{babel} 
\usepackage{t1enc}
\usepackage{fancyheadings} 
\usepackage{epic} 
\usepackage{eepic}

\newtheorem{prop}{Proposition} 
\newtheorem{lemma}[prop]{Lemma}

\newtheorem{thm}[prop]{Theorem}

\theoremstyle{definition}

\newcommand{\Ref}[1]{(\ref{#1})}
\newcommand{\findem}{\vspace{-.5cm} \begin{flushright} $\square~$ \end{flushright} \vspace{.4cm} }
\newcommand{\pare}[1]{\left( #1 \right)}
\newcommand{\parfrac}[2]{\left( \frac{#1}{#2} \right)}
\newcommand{\ite}{\noindent $\bullet~~$}
\newcommand{\iten}{\noindent -$~~$}

\newcommand{\Ss}{\textbf{S}}
\newcommand{\TT}{\textbf{T}}
\newcommand{\UU}{\textbf{U}}
\newcommand{\VV}{\textbf{V}}
\newcommand{\WW}{\textbf{W}}
\newcommand{\DD}{\textbf{D}}
\newcommand{\LL}{\textbf{L}}
\newcommand{\Fs}{\textbf{F}}
\newcommand{\GG}{\textbf{G}}
\newcommand{\HH}{\textbf{H}}
\newcommand{\KK}{\textbf{K}}

\newcommand{\EE}{\textbf{E}}

\newcommand{\SF}{\mathbb{S}}
\newcommand{\TF}{\mathbb{T}}
\newcommand{\UF}{\mathbb{U}}
\newcommand{\VF}{\mathbb{V}}
\newcommand{\WF}{\mathbb{W}}
\newcommand{\DF}{\mathbb{D}}
\newcommand{\LF}{\mathbb{L}}
\newcommand{\FF}{\mathbb{F}}
\newcommand{\GF}{\mathbb{G}}
\newcommand{\HF}{\mathbb{H}}
\newcommand{\KF}{\mathbb{K}}

\newcommand{\EF}{\mathbb{E}}
\newcommand{\XF}{\mathbb{X}}
\newcommand{\YF}{\mathbb{Y}}

\catcode`\@=11
\def\section{\@startsection{section}{1}%
 \z@{.7\linespacing\@plus\linespacing}{.5\linespacing}%
 {\normalfont\bfseries\scshape\centering}}

\def\subsection{\@startsection{subsection}{2}%
  \z@{.5\linespacing\@plus\linespacing}{.5\linespacing}%
  {\normalfont\bfseries\scshape}}

\def\subsubsection{\@startsection{subsubsection}{3}%
  \z@{.5\linespacing\@plus.7\linespacing}{-.5em}%
  {\normalfont\itshape}}
\catcode`\@=12
\title{On triangulations with high vertex degree}
\author{Olivier Bernardi}

\address{LaBRI, Universit\'e Bordeaux 1,
351 cours de la Lib\'eration,
  33405 Talence Cedex, France}
\email{bernardi@labri.fr}
\date{\today\\ Contact: \texttt{bernardi@labri.fr}}

\begin{document}
\begin{abstract}
We solve three enumerative problems concerning families of planar
maps. More precisely, we establish algebraic equations for the
generating function of  non-separable triangulations in which  all vertices have  
degree at least $d$, for a certain value $d$ chosen in $\{3,~4,~5\}$.\\
The originality of the problem lies in the fact that degree
restrictions are placed both on vertices and faces. Our proofs first follow Tutte's classical approach: 
we decompose maps by deleting the root and translate the decomposition into an equation satisfied by the generating function of the maps under consideration. 
Then we proceed to solve the equation obtained using a recent technique that extends the so-called \emph{quadratic method}. \\

\end{abstract}

\maketitle 
\thispagestyle{empty}

\section{Introduction}
The enumeration of planar maps (or \emph{maps} for short) has received a lot of attention in the combinatorists community for nearly fifty years. 
This field of research, launched by Tutte, was originally motivated by the four-color conjecture.
Tutte and his students considered a large number of map families corresponding to various constraints on face- or vertex-degrees. 
These seminal works, based on elementary decomposition techniques allied to a generating function approach, gave rise to many explicit results \cite{Tutte:census1,Tutte:census2,Tutte:census3,Tutte:census4,Mullin:triangulation-nonsep}. Fifteen years later, some physicists became interested in the subject and developed their own tools \cite{Bessis:matrix-integral,Brezin:matrix-integral,Hooft:matrix-integrals}  based on matrix integrals (see \cite{Zvonkin:intro} for an introduction). Their techniques proved very powerful for map enumeration \cite{DiFrancesco:census,DiFrancesco:gravity}.  More recently, a bijective approach based on \emph{conjugacy classes of trees} has emerged providing new insights on the subject \cite{Schaeffer:these,MBM:constellation,MBM:Ising2,Schaeffer:triangulation}.

However, when one considers a map family defined by both face- and vertex-constraints, each of the above mentioned methods seems relatively
ineffective and very few enumerative results are known.  There are  two major exceptions. 
First, and most importantly for this paper, the enumeration of \emph{triangulations} (faces have degree 3) in which all vertices have degree at least 3  and of triangulations \emph{without multiple edges} in which all vertices have degree at least 4 were performed  by Gao and Wormald using a compositional approach \cite{Gao:triang_degresup}. More recently, the enumeration of all \emph{bipartite} maps (faces have an even degree) according to the degree distribution of the vertices was accomplished using conjugacy classes of trees \cite{Schaeffer:eulerian,MBM:Ising2}. This result includes as a special case the enumeration of bipartite \emph{cubic} (vertices have degree 3) maps performed by Tutte via a generating function approach \cite{Tutte:census3,Tutte:chrom3}.\\

In this paper, we consider non-separable triangulations in which  all vertices have  degree at least $d$, for a certain value $d$ chosen in $\{3,~4,~5\}$. We establish algebraic equations for the generating function of each of these families. We also give the asymptotic behavior of the number of maps in each family. It is well-known that there is no triangulation in which  all vertices have  degree at least 6 (we shall prove this fact in Section \ref{section:preliminaires}). Hence, we have settled the problem of counting triangulations with 'high' vertex degree entirely. 

As mentioned above, the case $d=3$ was already solved  by Gao and Wormald \cite{Gao:triang_degresup}. Our proof differs from theirs.\\

Our proofs first follow Tutte's classical approach, which consists in translating the decomposition obtained by deletion of the root into a functional equation satisfied by the generating function. It is not clear at first sight why this approach should work here. As a matter of fact, finding a functional equation for triangulations with vertex degree at least 5 turns out to be rather complicated. But it eventually works if some of the constraints are relaxed at this stage of the solution. Our decomposition scheme requires to take into account, beside the size of the map, the degree of its root-face. Consequently, in order to write a functional equation, we need to consider a \emph{bivariate} generating function. We end up with an equation for the (bivariate) generating function in which the variable counting the degree of the root-face cannot be trivially eliminated. We then use a recent generalization of the \emph{quadratic method}  to get rid of the extra variable and compute an algebraic equation characterizing the \emph{univariate} generating function (see \cite{Brown:quadratic} and \cite[Section 2.9]{Goulden:Lagrange} for the quadratic method and \cite{MBM:quadratic} for its generalization).\\

This paper is organized as follows.  In Section \ref{section:preliminaires}, we recall some definitions on planar maps and introduce the main notations. In Section \ref{section:decomposition}, we recall the classic decomposition scheme  due to W.T. Tutte (by deletion of the root). We illustrate this scheme on the set of unconstrained non-separable near-triangulations. In Section \ref{section:functional_equations}, we apply the same decomposition scheme to the sets of near-triangulations in which any internal vertex has degree at least $3,~4,~5$. We obtain functional equation in which the variable $x$ counting the degree of the root-face cannot be trivially eliminated. In Section \ref{section:solution}, we use techniques generalizing the quadratic method in order to get rid of the variable $x$. We obtain algebraic equations for triangulations in which  \emph{any vertex not incident to the root} has degree at least $3,~4,~5$. In Section \ref{section:constraining}, we give algebraic equations for triangulations in which \emph{any vertex} has degree at least 3, 4. Lastly, in Section \ref{section:asymptotic} we study the asymptotic behavior of the number of maps in each family. \\

\section{Preliminaries and notations on maps} \label{section:preliminaires}
\begin{figure}[ht!]
\begin{center}
\input{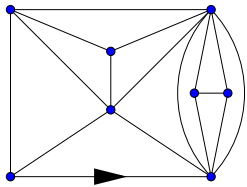} 
\caption{A non-separable near-triangulation.}  \label{fig:exp-triangulation}
\end{center}
\end{figure}

We begin with some vocabulary on maps.  A map is a proper embedding of a
connected graph into the two-dimensional sphere, considered up
to continuous deformations. A map is \emph{rooted} if one of its edges
is distinguished as the \emph{root} and attributed an
orientation. Unless otherwise specified, all maps under consideration
in this paper are rooted.  The face at the right of the root is called
the \emph{root-face} and the other faces are said
to be \emph{internal}. Similarly, the vertices incident to the root-face are
said to be \emph{external} and the others are said
to be \emph{internal}. Graphically, the root-face is usually represented as
the infinite face when the map is projected on the plane (see Figure
\ref{fig:exp-triangulation}). The endpoints of the root are
distinguished as its \emph{origin} and \emph{end} according to the
orientation of the root. A map is a \emph{triangulation} (resp. \emph{near-triangulation}) if all its faces (resp. all its internal faces) have degree 3.  For instance, the map of Figure \ref{fig:exp-triangulation} is a near-triangulation with root-face of degree 4.
Lastly, a map is \emph{non-separable} if it is loopless and  2-connected (the deletion of a vertex does not disconnect the map). For instance, the map in Figure \ref{fig:exp-triangulation} is non-separable.\\

In what follows, we enumerate 3 families of non-separable triangulations. We recall some basic facts about these maps. \vspace{.1cm}\\
\ite By definition, a non-separable triangulation has no loop. Therefore, the faces of non-separable triangulations are always homeomorphic to a triangle: they have three distinct vertices and three distinct edges. \vspace{.1cm}\\
\ite Consider a triangulation with $f$ faces, $e$ edges and $v$ vertices. Given the incidence relation between edges and faces, we have  $2e=3f$. Hence, the number of edges of a triangulation is a multiple of 3. Moreover, given the Euler relation ($v-e+f=2$), we see that a triangulation with $3n$ edges has $2n$ faces and $n+2$ vertices. \vspace{.1cm}\\
\ite Observe that a non-separable map (not reduced to an edge) cannot have a vertex of degree one. Let us now prove, as promised, that \emph{any triangulation has a vertex of degree less than 6}. Moreover, we prove that \emph{this vertex can be chosen not to be incident to the root}.  Indeed, if all vertices not incident to the root have degree  at least 6 the incidence relation between vertices and edges gives $2e \geq 6(v-2)+2$.  This contradicts the fact that triangulations with $e=3n$ edges have $v=n+2$ vertices.  This property shows that, if one considers the sets of non-separable triangulations with vertex degree at least $d$, \emph{the only interesting values of $d$ are $d=2$} (which corresponds to unconstrained non-separable triangulations) \emph{and $d=3,~4,~5$}. \\

\begin{figure}[ht]
\begin{center}
\input{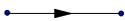} 
\caption{The link-map $L$.}  \label{fig:link-map}
\end{center}
\end{figure}

Let $\Ss$ be the set of non-separable rooted \emph{near-triangulations}. By convention, we exclude the map reduced to a vertex from $\Ss$. Thus, the smallest map in $\Ss$ is the map reduced to a straight edge (see Figure \ref{fig:link-map}). This map is called the \emph{link-map} and is denoted $L$. The vertices of other maps in $\Ss$ have degree at least 2. We consider three sub-families $\TT,~\UU,~\VV$ of $\Ss$.  The set $\TT$ (resp. $\UU$, $\VV$) is the subset of non-separable near-triangulations in which any \emph{internal} vertex  has degree at least 3 (resp. 4, 5). For each of the families $\WW=\Ss,~\TT,~\UU,~\VV$, we consider the bivariate generating function $\WF(x,z)$, where $z$ counts the size (the number of edges) and $x$ the degree of the root-face minus 2. That is to say,  $\WF(x)\equiv \WF(x,z)=\sum_{n,d} a_{n,d} x^{d}z^{n}$ where $a_{n,d}$ is the number of maps in $\WW$ with size $n$ and root-face of degree $d+2$.
For instance, the link-map $L$, which is the smallest map in all our families, has contribution $z$ to the generating function. Therefore, $\WF(x)=z+o(z)$. Since the degree of the root-face is bounded by two times the number of edges, the generating function $\WF(x,z)$ is a power series in the main variable $z$ with polynomial coefficients in the secondary variable $x$. For each family $\WW=\Ss,~\TT,~\UU,~\VV$, we will characterize the generating function $\WF(x)$ as the unique power series solution of a functional equation (see Equation \Ref{eq:Sintro} and Propositions  \ref{prop:T}, \ref{prop:U}, \ref{prop:V}).\\

We also consider the set $\Fs$ of non-separable rooted  \emph{triangulations} and three of its subsets $\GG~,\HH~,\KK$. The set $\GG$ (resp. $\HH$, $\KK$) is the subset of non-separable  triangulations in which any vertex \emph{not incident to the root} has degree at least 3 (resp. 4, 5). 
As observed above, the number of edges of a triangulation is always a multiple of 3. 
To each of the families $\LL=\Fs,~\GG,~\HH,~\KK$, we associate the univariate generating function $\LF(t)=\sum_{n} a_{n} t^{n}$ where $a_n$ is the number of maps in $\LL$ with $3n$ edges ($2n$ faces and $n+2$ vertices).  For each family we will 
give an algebraic equation satisfied by $\LF(t)$ (see Equation \Ref{eq:F} and Theorems  \ref{thm:G}, \ref{thm:H}, \ref{thm:K}).\\

\begin{figure}[ht!]
\begin{center}
\input{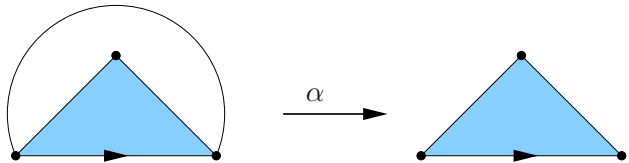}
\caption{Near-triangulations rooted on a digon and triangulations.} \label{fig:alphafig}
\end{center}
\end{figure}

There is a simple connection between the generating functions  $\FF(t)$ (resp. $\GF(t)$, $\HF(t)$, $\KF(t)$)
and $\SF(x)$ (resp. $\TF(x)$, $\UF(x)$, $\VF(x)$). 
Consider a non-separable near-triangulation distinct from $L$ rooted on a digon (i.e. the root-face has degree 2). Deleting  the external edge that is not the root produces a non-separable triangulation (see  Figure \ref{fig:alphafig}). This mapping establishes a one-to-one correspondence between the set of triangulations $\Fs$ (resp. $\GG~,\HH~,\KK$) and the set of near-triangulations in $\Ss-\{L\}$ (resp. $\TT-\{L\},~\UU-\{L\},~\VV-\{L\}$) rooted on a digon.\\
\vspace{-.2cm}  \vspace{-.3cm} 

\noindent For $\WW \in \{\Ss,~\TT,~\UU,~\VV\}$, the power series $\WF(0)\equiv \WF(0,z)$ is the generating function of near-triangulations in $\WW$ rooted on a digon.  Given that the link-map has contribution $z$, we have  
\begin{eqnarray} \SF(0)=z+z\FF(z^3), \hspace{.2cm} \TF(0)=z+z\GF(z^3),\hspace{.2cm}  \UF(0)=z+z\HF(z^3),\hspace{.2cm}  \VF(0)=z+z\KF(z^3). \label{eq:alpha} \end{eqnarray}\\

\section{The decomposition scheme}\label{section:decomposition}
In the following, we adopt Tutte's classical approach for enumerating maps. That is, we decompose maps by deleting their root and translate this combinatorial decomposition into an equation satisfied by the corresponding generating function. 
In this section we illustrate  this approach on  unconstrained non-separable triangulations (this was first done in \cite{Mullin:triangulation-nonsep}). We give all the details on this simple case in order to prepare the reader to the more complicated cases of constrained non-separable triangulations treated in the next section.

\begin{figure}[ht!]
\begin{center}
\input{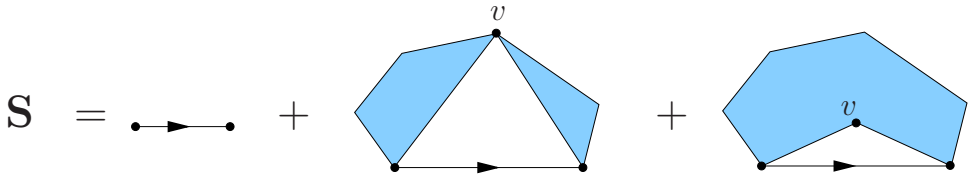} 
\caption{Decomposition of non-separable near-triangulations.} \label{fig:decompose}
\end{center}
\end{figure}

We recall that $\Ss$ denotes the set of non-separable near-triangulations and $\SF(x)=\SF(x,z)$ the corresponding generation function. As observed before, the link-map $L$ has contribution $z$ to the generating function $\SF(x)$. We decompose the other maps by deleting the root. Let $M$ be a non-separable triangulation distinct from $L$. Since $M$ is non-separable, the root of $M$ is not an isthmus. Therefore, the face at the left of the root is internal, hence has degree 3. Since $M$ has no loop, the three vertices incident to this face are distinct. We denote by  $v$ the vertex not incident to the root. When analyzing what can happen to $M$ when deleting its root, one is led to distinguish two cases (see Figure \ref{fig:decompose}). \\
\vspace{-.2cm}  \vspace{-.3cm} 

\begin{figure}[ht!]
\begin{center}
\input{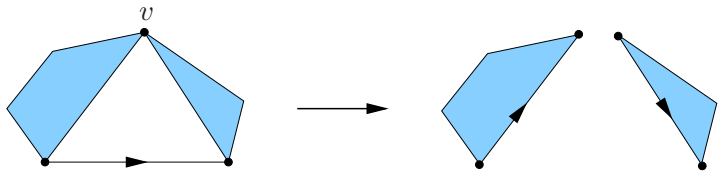}
\caption{Case 1. The vertex $v$ is incident to the root-face.} \label{fig:decompose-cas1}
\end{center}
\vspace{-.5cm} 
\end{figure}

\begin{figure}[ht!]
\begin{center}
\input{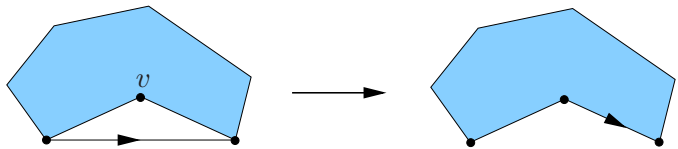}
\caption{Case 2. The vertex $v$ is not incident to the root-face.}
\label{fig:decompose-cas2}
\end{center}
\end{figure}

\noindent Either the vertex $v$ is incident to the root-face, in which case the map obtained by deletion of the root is separable (see Figure \ref{fig:decompose-cas1}). Or $v$ is not incident to the root-face and the map obtained by deletion of the root is a non-separable near-triangulation (see Figure \ref{fig:decompose-cas2}). In the first case, the map obtained is in correspondence with an ordered pair of non-separable near-triangulations. This correspondence is bijective, that is, any ordered pair is the image of exactly one near-triangulation.  In the second case the degree of the root-face is increased by one. Hence the root-face of the near-triangulation obtained has degree at least 3. Here again, any near-triangulation in which the root-face has degree at least 3 is the image of exactly one near-triangulation.\\
  
We want to translate this analysis into a functional equation. Observe
that the degree of the root-face appears in this analysis. This is why
we are \emph{forced} to introduce the variable $x$ counting this
parameter in our generating function $\SF(x,z)$.  For this reason, following  Zeilberger's terminology \cite{Zeilberger:catalytic}, 
the secondary variable $x$ is said to be \emph{catalytic}: we need it to write the functional equation, but we shall try to get rid of it later. \\

In our case, the decomposition easily translates into the following equation (details will be given in Section \ref{section:functional_equations}):
\begin{eqnarray} \SF(x,z)=z+x z\SF(x,z)^2+\frac{z}{x}\left( \SF(x,z)-\SF(0,z)\right). \label{eq:Sintro} \end{eqnarray}
The first summand of the right-hand side accounts for the link map, the second summand corresponds to the case in which the vertex $v$ is incident to the root-face, and the third summand corresponds to the case in which $v$ is not incident to the root-face.\\

It is an easy exercise to check that this equation defines the series $\SF(x,z)$ uniquely as a power series in $z$ with polynomial coefficients in $x$. By solutions techniques  presented in Section  \ref{section:solution}, we can derive from Equation  \Ref{eq:Sintro} a polynomial equation  satisfied by the series $\SF(0,z)$ where the extra variable $x$ does not appear anymore. This equation reads
\begin{eqnarray} \SF(0,z)~=~z-27z^4+36z^3\SF(0,z)-8z^2\SF(0,z)^2-16z^4\SF(0,z)^3.\label{eq:S(0)}\end{eqnarray} 
Given that $\SF(0,z)=z+z\FF(z^3)$, we deduce the algebraic equation
\begin{eqnarray}  \FF(t)=t(1-16t)-t(48t-20)\FF(t)-8t(6t+1)\FF(t)^{2}-16t^{2}\FF(t)^{3},\label{eq:F}\end{eqnarray} 
characterizing $\FF(t)$ (the generating function of non-separable triangulations) uniquely as a power series in $t$. From this equation one can derive the asymptotic behavior of the coefficients of $\FF(t)$, that is, the number of non-separable triangulations of a given size (see Section \ref{section:asymptotic}). \\

\section{Functional equations}\label{section:functional_equations}

\begin{figure}[ht!]
\begin{center}
\input{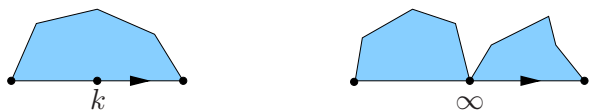}
\caption{Generic elements of the sets  $\WW_{k}$ and $\WW_{\infty}$.} \label{fig:S-k-infini}
\end{center}
\end{figure}

In this section, we apply the decomposition scheme presented in Section \ref{section:decomposition} to the families $\TT,~\UU,~\VV$ of non-separable near-triangulations in which all internal vertices have degree at least 3,~4,~5. We obtain functional equations 
 satisfied by the corresponding generating functions $\TF(x),~\UF(x),~\VF(x)$.\\
Note that, when one deletes the root of a map, the degree of its endpoints is lowered by one. Given the decomposition scheme, this remark explains why we are led to consider the near-triangulations where only \emph{internal} vertices have a degree constraint.
However, we need to control the degree of the root's origin since it may come from an internal vertex (see Figure \ref{fig:decompose-cas2}). This leads to the following notations. 
Let $\WW$ be one of the sets $\Ss,~\TT,~\UU,~\VV$. We define $\WW_k$ as the set of maps in $\WW$ such that the root-face has degree at least 3 and the origin of the root has degree $k$. We also define the set $\WW_\infty$ as the set of (separable) maps obtained by gluing the root's end  of a map in $\WW$ with the root's origin of a map in $\WW$.  The root of the map obtained is chosen to be the root of the second map.
Generic elements of the sets $\WW_{k}$ and $\WW_{\infty}$ are shown in Figure \ref{fig:S-k-infini}. We also write $\displaystyle \WW_{\geq k}\triangleq \WW_\infty \cup \bigcup_{j\geq k} \WW_j$. The notation $\WW_{\geq k}$, which at first sight might seem awkward, allows to unify the two possible cases of our decomposition scheme (Figure \ref{fig:decompose-cas1} and \ref{fig:decompose-cas2}). It shall simplify our arguments and equations (see for instance Equations (\ref{eq:decomposeS}-\ref{eq:decomposeV})) which will prove a valuable property.  \\

The symbols $\WF_k(x,z)$, $\WF_\infty(x,z)$  and $\WF_{\geq k}(x,z)$ denote the bivariate generating functions of the sets $\WW_k$, $\WW_\infty$  and $\WW_{\geq k}$ respectively. In these series, as in $\WF(x,z)$, the contribution of a map with $n$ edges and root-face degree $d+2$ is $x^{d}z^n$.\\

\begin{figure}[ht!]
\begin{center}
\input{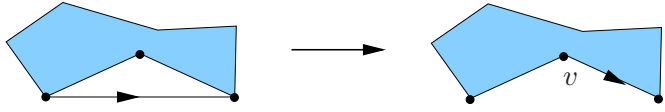}
\caption{Mapping induced by deletion of the root: the  vertex $v$ can be a separating point in which case the map is in $\WW_\infty$.}
\label{fig:decompose-cas-double}
\end{center}
\end{figure}

We are now ready to apply the decomposition scheme to the triangulations in $\TT,~\UU,~\VV$. Consider a near-triangulation $M$ distinct from $L$ in $\WW=\Ss,~\TT,~\UU,~\VV$. As observed before, the face at the left of the root is an internal face incident to three distinct vertices. We denote by $v$ the vertex not incident to the root. 
If $v$ is external, the  deletion of the root produces a map in $\WW_\infty$ (see Figure \ref{fig:decompose-cas-double}). If $v$ is internal and $M$ is in $\Ss$ (resp. $\TT,~\UU,~\VV$) then $v$ has degree at least $2$ (resp. $3,~4,~5$) and the  map obtained by deleting the root is in $\bigcup_{k \geq 2} \Ss_k$ (resp. $\bigcup_{k \geq 3} \TT_k$, $\bigcup_{k \geq 4} \UU_k$, $\bigcup_{k \geq 5} \VV_k$). Therefore, the deletion of the root induces a mapping  from $\Ss-\{L\}$ (resp. $\TT-\{L\}$, $\UU-\{L\}$, $\VV-\{L\}$) to $\Ss_{\geq 2}$ (resp.  $\TT_{\geq 3}$, $\UU_{\geq 4}$, $\VV_{\geq 5}$).\\ 

\noindent This mapping is clearly bijective. Moreover, the map obtained after deleting the root has size lowered by one and root-face degree increased by one. This analysis translates into the following equations:
\begin{eqnarray}
\SF(x)&=&z+\frac{z}{x}\SF_{\geq 2}(x)~, \label{eq:decomposeS}  \\ 
\TF(x)&=&z+\frac{z}{x}\TF_{\geq 3}(x)~,  \label{eq:decomposeT} \\
\UF(x)&=&z+\frac{z}{x}\UF_{\geq 4}(x)~,  \label{eq:decomposeU} \\
\VF(x)&=&z+\frac{z}{x}\VF_{\geq 5}(x)~.  \label{eq:decomposeV}
\end{eqnarray}

In view of Equation (\ref{eq:decomposeS}), we will obtain 
a non-trivial equation for $\SF(x)$ if we can express $\SF_{\geq
  2}(x)$ in terms of $\SF(x)$. Similarly, we will obtain a
non-trivial equation for $\TF(x)$ if we can express $\TF_{\geq 2}(x)$
and $\TF_{2}(x)$ in terms of $\TF(x)$. Similar statements hold for $\UF(x)$ and $\VF(x)$.
Thus, our first task will be to evaluate  $\WF_{\geq 2}(x)$ for $\WF$ in $\{\SF,~\TF,~\UF,~\VF\}$.\\

By definition, $\WW_\infty$ is in bijection with $\WW^2$, which translates into the functional equation
$$ \WF_\infty(x) ~=~ x^2\WF(x)^2~.$$ 
Observe that $ \bigcup_{k\geq 2} \WW_k$ is the set of maps in $\WW$ for which the root-face has degree at least 3, that is, all maps except those rooted on a digon. Since  $\WF(0)$ is the generating function of maps in $\WW$ rooted on a digon, we have 
$$ \sum_{k \geq 2} \WF_k(x) ~=~ \WF(x)-\WF(0). $$  
Given that $ \WW_{\geq 2}=\WW_\infty \cup  \bigcup_{k \geq 2} \WW_k ~,$ we obtain, for $\WF$ in $\{\SF,~\TF,~\UF,~\VF\}$,
\begin{eqnarray} \WF_{\geq 2}(x)~=~x^2\WF(x)^2 ~+~ \left(\WF(x)-\WF(0) \right) \hspace{1.5cm} \textrm{for }\WF ~\textrm{in}~\{\SF,~\TF,~\UF,~\VF\}~.  \label{eq:Rsup2toR} \end{eqnarray} 

Equations \Ref{eq:decomposeS} and \Ref{eq:Rsup2toR} already prove Equation \Ref{eq:Sintro} announced in Section \ref{section:decomposition}:
$$ \SF(x)=z+x z\SF(x)^2+z\parfrac{\SF(x)-\SF(0)}{x}~.$$

In order to go further, we need to express $\TF_2(x),~\UF_2(x),~\UF_3(x),~\VF_2(x),~\VF_3(x)$ and $\VF_4(x)$ (see Equations (\ref{eq:decomposeT}-\ref{eq:decomposeV})).  We begin with an expression of  $\WF_2(x)$ for  $\WF$ in $\{\SF,~\TF,~\UF,~\VF\}$. 
Observe that for $\WW=\{\Ss,~\TT,~\UU,~\VV\}$, the set  $\WW_2$ is in bijection with $\WW$ by the mapping illustrated in Figure \ref{fig:R2toR}. Consequently we can write 
\begin{eqnarray} \WF_{2}(x)~=~xz^2\WF(x)~ \hspace{1.5cm} \textrm{for }\WF ~\textrm{in}~\{\SF,~\TF,~\UF,~\VF\} ~. \label{eq:R2toR} \end{eqnarray}

\vspace{-.2cm} \begin{figure}[ht!]
\begin{center}
\input{R2toR.pstex_t}
\caption{A bijection between  $\WW_2$ and $\WW$.} \label{fig:R2toR}
\end{center}
\end{figure} \vspace{-.3cm}

This suffices to obtain an equation for the set $\TT$:
$$\hspace{-.2cm}
\begin{array}{lllr}
\TF(x)&=&\displaystyle z+\frac{z}{x}\TF_{\geq 3}(x) & \hspace{.1cm} \textrm{by } \Ref{eq:decomposeT} \vspace{.1cm}\\ 
&=&\displaystyle z+\frac{z}{x}\left(\TF_{\geq 2}(x) - \TF_{2}(x)\right) \vspace{.1cm} \\ 
&=&\displaystyle z+\frac{z}{x} \left(x^2\TF(x)^2 + \left(\TF(x)-\TF(0) \right) - xz^2\TF(x)\right) & \hspace{.2cm} \textrm{by } \Ref{eq:Rsup2toR} \textrm{ and }  \Ref{eq:R2toR}.
\end{array}$$

\begin{prop} \label{prop:T} 
The generating function $\TF(x)$ of non-separable near-triangulations in which all internal vertices have degree at least 3 satisfies:
\begin{eqnarray}\TF(x)&=&z~+~ xz\TF(x)^2 ~+~ z\left(\frac{\TF(x)-\TF(0)}{x} \right) ~-~ z^3\TF(x)~. \label{eq:T} \end{eqnarray}
\end{prop}

In order to  find an equation concerning the sets $\UU$ and $\VV$, we now need to express $\UF_3(x)$ and $\VF_3(x)$ in terms of $\UF(x)$ and $\VF(x)$ respectively. Let $\WW$ be $\UU$ or $\VV$ and $M$ be a map in $\WW_3$. By definition, the root-face of $M$ has degree at least 3 and its root's origin $u$ has degree 3. We denote by $a$ and $b$ the vertices preceding and following $u$ on the root-face (see Figure \ref{fig:U3toUsup3}). Since the map $M$ is non-separable, the vertices $a, b$ and $u$ are distinct. Let $v$ be the third vertex adjacent to $u$. Since $M$ cannot have loops, $v$ is distinct from $a, b, v$.\\
Suppose that $M$ is in $\UU_3$ (resp. $\VV_3$)  and consider the operation of deleting $u$ and the three adjacent edges. If the vertex $v$ is internal it has degree $d\geq 4$ (resp.  $d\geq 5$) and the map obtained is in $\UU_{d-1}$ (resp. $\VV_{d-1}$). If it is external, the map obtained is in $\UU_\infty$ (resp. $\VV_\infty$). Thus, the map obtained is in $\UU_{\geq 3}$ (resp. $\VV_{\geq 4}$). This correspondence is clearly bijective. It gives 
\begin{eqnarray} 
\UF_3(x)&=&z^3\UF_{\geq 3}(x)~=~z^3(\UF_{\geq 2}(x) - \UF_{2}(x))~,\label{eq:U3} \\
\VF_3(x)&=&z^3\VF_{\geq 4}(x)~=~z^3(\VF_{\geq 2}(x) - \VF_{2}(x) - \VF_3(x))~.\label{eq:V3bis}
\end{eqnarray}

\vspace{-.2cm} \begin{figure}[ht!]
\begin{center}
\input{U3toUsup3.pstex_t}
\caption{A bijection between $\UU_3$ and $\UU_{\geq 3}$ (resp. $\VV_3$ and $\VV_{\geq 4}~$).} \label{fig:U3toUsup3}
\end{center}
\end{figure} \vspace{-.3cm} 

We are now ready to establish the functional equation concerning $\UU$:
$$\hspace{-.2cm}
\begin{array}{lllr}
\UF(x)&=&\displaystyle z+\frac{z}{x}\UF_{\geq 4}(x) & \hspace{.2cm} \textrm{by } \Ref{eq:decomposeU} \vspace{.1cm}\\ 
&=&\displaystyle z+\frac{z}{x}\left(\UF_{\geq 2}(x) - \UF_{2}(x)  - \UF_{3}(x) \right)  \vspace{.1cm}\\ 
&=&\displaystyle z+\frac{z(1-z^3)}{x}\left(\UF_{\geq 2}(x)  - \UF_{2}(x)\right)  & \hspace{.1cm} \textrm{by } \Ref{eq:U3} \vspace{.1cm}\\ 
&=&\displaystyle z+\frac{z(1-z^3)}{x} \left(x^2\UF(x)^2 + \left(\UF(x)-\UF(0) \right) - xz^2\UF(x)\right) & \hspace{.1cm} \textrm{by } \Ref{eq:Rsup2toR} \textrm{ and } \Ref{eq:R2toR}.
\end{array}$$ 

\begin{prop} \label{prop:U} 
The generating function $\UF(x)$ of non-separable near-triangulations in which all internal vertices have degree at least 4 satisfies:
\begin{eqnarray}\hspace{-.2cm} \UF(x)=z+ xz(1-z^3)\UF(x)^2 + z(1-z^3)\left(\frac{\UF(x)-\UF(0)}{x} \right) - z^3(1-z^3)\UF(x)~. \label{eq:U} \end{eqnarray}
\end{prop}


We proceed to find an equation concerning the set $\VV$. This will require significantly more work than the previous cases. We write
\begin{eqnarray}\VF(x)&=&z~+~\frac{z}{x}\VF_{\geq 4}(x)~=~ z~+~\frac{z}{x}\left( \VF_{\geq 2}(x)- \VF_{2}(x)- \VF_{3}(x)- \VF_{4}(x)\right) \label{eq:redecomposeV} \end{eqnarray}
and we want to express $\VF_{\geq 2}(x),~ \VF_{2}(x),~ \VF_{3}(x)$ and $\VF_{4}(x)$ in terms of $\VF(x)$.
We already have such expressions for   $\VF_{\geq 2}(x)$ and $\VF_{2}(x)$ (by Equations \Ref{eq:Rsup2toR} and \Ref{eq:R2toR}). Moreover, Equation \Ref{eq:V3bis} can be rewritten as 
\begin{eqnarray} \VF_3(x)&=&\frac{z^3}{1+z^3}\left(\VF_{\geq 2}(x) - \VF_{2}(x) \right)~.\label{eq:V3} \end{eqnarray}

It remains to express $\VF_{4}(x)$ in terms of $\VF(x)$. Unfortunately, this requires some efforts and some extra notations.
We define $\VV_{k,l}$ as the set of maps in $\VV$ such that the root-face has degree at least 4, the root's origin has degree $k$ and the root's end has degree $l$ (see  Figure \ref{fig:S-infini-infini}). 
The set  $\VV_{k,\infty}$ is the set of maps obtained by gluing  the root's end  of a map in $\VV_k$ with the root's origin of a map in $\VV$. The root of the new map obtained is the root of the map in $\VV_k$. 
The set  $\VV_{\infty,k}$  is the set of maps obtained by gluing the root's end of a map in $\VV$ with the root's origin of map in $\VV$ for which the root-face has degree at least 3 and the root's end has degree $k$. The root of the new map obtained is the root of the second map.  The set  $\VV_{\infty,\infty}$ is obtained by gluing 3 maps of $\VV$ as indicated in Figure \ref{fig:S-infini-infini}.\\
\begin{figure}[ht!]
\begin{center}
\input{S-infini-infini.pstex_t}
\caption{The sets  $\VV_{k,l}$, $\VV_{\infty,k}$, $\VV_{k,\infty}$ and $\VV_{\infty,\infty}$.} \label{fig:S-infini-infini}
\end{center}
\end{figure} 

\noindent We also write  $\displaystyle \VV_{k, \geq l}\triangleq \bigcup_{i\geq l} \VV_{k,i} \cup \VV_{k,\infty}$ and  
$$\displaystyle \VV_{\geq k, \geq l}\triangleq \bigcup_{i\geq k,j\geq l} \VV_{i,j} \cup \bigcup_{i\geq k} \VV_{i,\infty}  \cup \bigcup_{j\geq l}\VV_{\infty,j}  \cup \VV_{\infty,\infty}.$$ 
As before, if $\WW$ is any of these sets, the symbol $\WF$ denotes the corresponding generating function, where the contribution of a map of size $n$ and root-face degree $d+2$ is $x^{d}z^n$.\\

Moreover, we consider the subset $\DD$ of $\VV$ composed of maps for which the root-face is a digon. The set of maps in $\DD$ for which the root-vertex has degree $k$ will be denoted by $\DD_k$. We write $\DD_{\geq k}=\bigcup_{j\geq k} \DD_j$. Lastly, if $\EE$ is one of the set $\DD$, $\DD_k$ or $\DD_{\geq k}$, the symbol $\EF$ denotes the corresponding (univariate) generating function, where the contribution of a map of size $n$ is $z^n$. As observed before,  $\DF=\VF(0)$. \\

\begin{figure}[ht!]
\begin{center}
\input{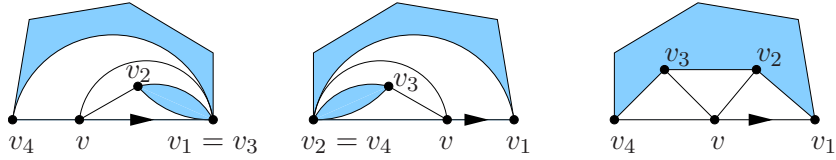}
\caption{Three configurations for a map in $\VV_4$.} \label{fig:decomposeV4}
\end{center}
\end{figure}

We can now embark on the decomposition of $\VV_4$. We consider a map $M$ in $\VV_4$ with root-vertex $v$. By definition, $v$ has degree 4. Let $e_1,~e_2,~e_3,~e_4$ be the edges incident to $v$  in counterclockwise order  starting from the root $e_1$. We denote by $v_i,~i=1\ldots4$ the endpoint of $e_i$ distinct from $v$. Since $M$ is non-separable and its root-face has degree at least $3$, the vertices $v_1$ and $v_4$ are distinct. Moreover, since $M$ has no loop we have $v_1 \neq v_2$,  $v_2 \neq v_3$ and $v_3 \neq v_4$.
Therefore, only three configurations are possible: either $v_1=v_3$, the two other vertices being distinct, or symmetrically, $v_2=v_4$, the other vertices being distinct, or $v_1,v_2,v_3,v_4$ are all distinct. The three cases are illustrated in Figure \ref{fig:decomposeV4}.
 
\begin{figure}[ht!]
\begin{center}
\input{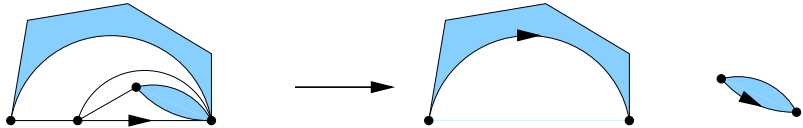}
\caption{A bijection between maps of the first type in $\VV_4$ and $\VV\times \DD_{\geq 4}$.} \label{fig:decomposeV4-cas1}
\end{center}
\end{figure}

\begin{figure}[ht!]
\begin{center}
\input{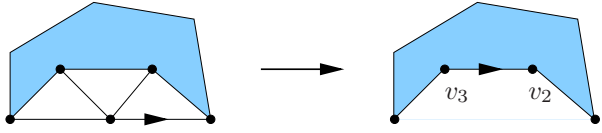}
\caption{A bijection between maps of the third type in $\VV_4$ and $\VV_{\geq 4, \geq 4}$.} \label{fig:decomposeV4-cas2}
\end{center}
\end{figure}

In the case $v_1=v_3$, the map can be decomposed into an ordered pair of maps in $\VV\times \DD_{\geq 4}$ (see Figure \ref{fig:decomposeV4-cas1}). This decomposition is clearly bijective. The symmetric case $v_2=v_4$  admits a similar treatment. In the last case ($v_1,v_2,v_3,v_4$ all distinct) the map obtained  from $M$ by deleting $e_1,~e_2,~e_3,~e_4$ is in $\VV_{\geq 4, \geq 4}$ (see Figure \ref{fig:decomposeV4-cas2}). Note that this case contains several subcases depending on whether $v_2$ and $v_3$ are separating points or not. But again the correspondence is clearly bijective.\\

\noindent This correspondence gives 
\begin{eqnarray} \VF_4(x)~=~2xz^4\VF(x)\DF_{\geq 4}~ + ~ \frac{z^4}{x} \VF_{\geq 4,\geq 4}(x)~. \label{eq:decomposeV4} \end{eqnarray}

\begin{figure}[ht!]
\begin{center}
\input{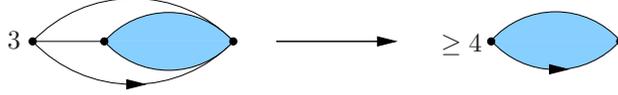}
\caption{A bijection between $\DD_3$ and $\DD_{\geq 4}$.} \label{fig:D3toDsup4}
\end{center}
\end{figure}

It remains to express the generating functions $\DF_{\geq 4}$ and $\VF_{\geq 4,\geq 4}(x)$ in terms of $\VF(x)$. We start with $\DF_{\geq 4}$.\\
We have $\DF_{\geq 4}=\DF - \DF_1 - \DF_2 - \DF_3~$. We know that $\DF=\VF(0)$. Moreover, the set $\DD_1$ only contains the link-map and $\DD_2$ is empty. Hence $\DF_1=z$ and $\DF_2=0$. Lastly, the set $\DD_3$ is in correspondence with $\DD_{\geq 4}$ by the bijection represented in Figure \ref{fig:D3toDsup4}.  This gives $\DF_3=z^3\DF_{\geq 4}$.

\noindent Putting these results together and solving for $\DF_{\geq 4}$,  we get 
\begin{eqnarray} \DF_{\geq 4}&=&\frac{\VF(0)-z}{1+z^3}~.\label{eq:D4}\end{eqnarray}

We now want to express the generating function $\VF_{\geq 4, \geq 4}(x)$. We first divide our problem as follows (the equation uses the trivial bijections between the sets $\VV_{\alpha,\beta}$ and $\VV_{\beta,\alpha}$) :
\begin{eqnarray} \VF_{\geq 4, \geq 4}(x)=  \VF_{\geq 2, \geq 2}(x) -\VF_{2, 2}(x)-2 \VF_{2, \geq 3}(x)-\VF_{3, 3}(x) - 2 \VF_{3, \geq 4}(x)~.\label{eq:decomposeUsup4sup4} \end{eqnarray}
We now treat separately the different summands in the right-hand-side of this equation.

\ite $\VV_{\geq 2, \geq 2}$ :
It follows easily from the definitions that :
$$\displaystyle \VF_{\geq 2, \geq 2}(x)=\sum_{k\geq 2, l\geq 2} \VF_{k,l}(x)+ 2 \sum_{k\geq 2} \VF_{\infty,k}(x) + \VF_{\infty,\infty}(x).$$
 \iten The set $\bigcup_{k\geq 2, l\geq 2} \VV_{k,l}$ is the set of maps in $\VV$ for which the root-face has degree at least 4. Thus, 
$$\sum_{k \geq 2,l\geq2}\VF_{k,l}(x)~=~\VF(x)-\VF(0)-x[x]\VF(x),$$
where $[x]\VF(x)$ is the coefficient of $x$ in $\VF(x)$.\\
\iten By definition, the set $\bigcup_{k\geq 2} \VV_{\infty,k}$ is in bijection with $\VV\times \bigcup_{k\geq 2}\VV_k$. Moreover, the set  $\bigcup_{k\geq 2}\VV_k$ is the set of maps in $\VV$ for which the root-face has degree at least 3. This gives
$$\sum_{k \geq 2}\VF_{\infty,k}(x)~=~x^2\VF(x)\left(\VF(x)-\VF(0)\right).$$
\iten By definition, the set $\VV_{\infty,\infty}$ is in bijection with $\VV^3$, which gives 
$$\VF_{\infty,\infty}(x)~=~x^4\VF(x)^3.$$
 
\noindent Summing these contributions we get 
\begin{eqnarray}\VF_{\geq 2, \geq 2}(x)= \VF(x)-\VF(0)-x[x]\VF(x)+2x^2\VF(x)\left(\VF(x)-\VF(0)\right)+x^4\VF(x)^3. \label{eq:Vsup2sup2}\end{eqnarray}

\ite $\VV_{2, 2}$ :
The set $\VV_{2, 2}$ is empty (the face at the left of the root would be of degree at least 4), hence 
\begin{eqnarray}\VF_{2, 2}(x)~=~ 0~. \label{eq:V22}\end{eqnarray}

\begin{figure}[ht!]
\begin{center}
\input{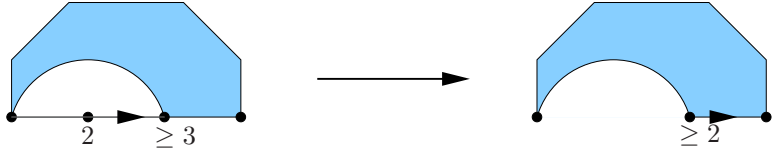}
\caption{A bijection between $\VV_{2, \geq 3}$ and $\VV_{\geq 2}$.} \label{fig:V2sup3toVsup2}
\end{center}
\end{figure}

\ite $\VV_{2, \geq 3}$ :
The set $\VV_{2, \geq 3}$  is in bijection with $\VV_{\geq 2}$ by the mapping illustrated in Figure \ref{fig:V2sup3toVsup2}. This gives $\VF_{2, \geq 3}(x)~=~ xz^2\VF_{\geq2}(x)$. From this, using  Equation \Ref{eq:Rsup2toR}, we obtain 
\begin{eqnarray}\VF_{2, \geq 3}(x)~=~  xz^2(\VF(x)-\VF(0)+x^2\VF(x)^2). \label{eq:V2sup3}\end{eqnarray}

\begin{figure}[ht!]
\begin{center}
\input{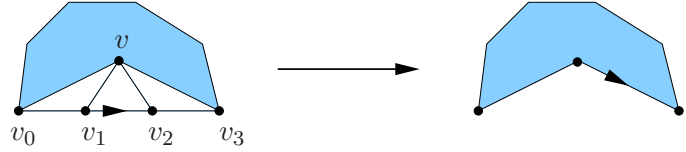}
\caption{A bijection between $\VV_{3, 3}$ and $\VV_{\geq 3}$.} \label{fig:V33toVsup3}
\end{center}
\end{figure}

\ite $\VV_{3,3}$ :
We consider a map  $M$ in $\VV_{3,3}$. We denote by $v_1$ the root's origin, $v_2$ the root's end,  $v_0$ the vertex preceding $v_1$ on the root-face and $v_3$ the vertex following $v_2$ (see Figure \ref{fig:V33toVsup3}). Since $M$ is non-separable and its root-face has degree at least 4, the vertices $v_i,~i=1\ldots4$ are all distinct. The third vertex $v$ adjacent with $v_1$ is also the third vertex adjacent  with $v_2$ (or the face at the left of the root would not be a triangle). Since $M$ has no loop, $v$ is distinct from  $v_i,~i=1\ldots 4$. 
From these considerations, it is easily seen that the set $\VV_{3,3}$ is in bijection with the set $\VV_{\geq 3}$ by the mapping illustrated in Figure \ref{fig:V33toVsup3}. (Note that this correspondence includes two subcases depending on $v$ becoming a separating point of not). We obtain 
$$\VF_{3,3}(x)~=~ xz^5\VF_{\geq3}(x)~=~xz^5\left(\VF_{\geq2}(x)-\VF_{2}(x)\right).$$
From this, using Equations \Ref{eq:Rsup2toR} and \Ref{eq:R2toR}, we get
\begin{eqnarray}\VF_{3, 3}(x)~=~  xz^5(\VF(x)-\VF(0)+x^2\VF(x)^2-xz^2\VF(x)). \label{eq:V33}\end{eqnarray}

\begin{figure}[ht!]
\begin{center}
\input{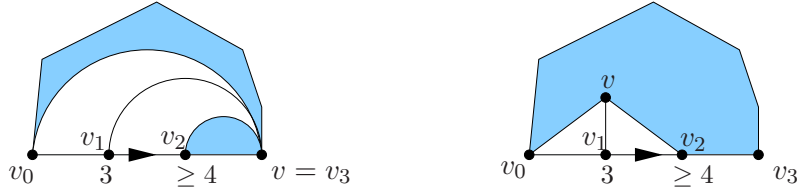}
\caption{Two configurations for a map in $\VV_{3,\geq 4}$.} \label{fig:decomposeV3sup4}
\end{center}
\end{figure}

\ite $\VV_{3,\geq 4}$ : 
Let $M$ be a map in  $\VV_{3,\geq 4}$. We denote by $v_1$ the root's origin, $v_2$ the root's end,  $v_0$ the vertex preceding $v_1$ on the root-face and $v_3$ the vertex following $v_2$ (see Figure \ref{fig:decomposeV3sup4}). Since $M$ is non-separable and its root-face has degree at least 4, the vertices $v_i,~i=1\ldots4$ are all distinct. Let $v$ be the third vertex adjacent to $v_1$. 
Two cases can occur: either  $v=v_3$ in which case the map decomposes into an ordered pair of maps in $\VV\times \DD_{\geq 3}$, or $v$ is distinct from $v_i,~i=1\ldots 4$ in which case the map is in correspondence with a map in $\VV_{\geq 4,\geq 3}$ (this includes two subcases  depending on $v$ becoming a separating point of not).
In both cases the correspondence is clearly bijective. This  gives
$$\VF_{3, \geq 4}(x)~=~  x^2z^3\VF(x)\DF_{\geq 3} ~+~z^3\VF_{\geq 4, \geq 3}(x).$$ 

\noindent Given that $\DF_{\geq 3}~=~\DF -\DF_1-\DF_2~=~\VF(0)-z$, we obtain
$$\VF_{3, \geq 4}(x)~=~  x^2z^3\VF(x)(\VF(0)-z) ~+~z^3(\VF_{\geq 4, \geq 4}(x)+ \VF_{3, \geq 4}(x)),$$
and solving for $\VF_{3, \geq 4}(x)$ we get 
\begin{eqnarray}\VF_{3, \geq 4}(x)~=~  \frac{z^3}{1-z^3}\left(x^2\VF(x)(\VF(0)-z) ~+~\VF_{\geq 4, \geq 4}(x)\right). \label{eq:V3sup4}\end{eqnarray}

We report Equations (\ref{eq:Vsup2sup2} - \ref{eq:V3sup4}) in Equation \Ref{eq:decomposeUsup4sup4} and solve for $\VF_{\geq 4,\geq 4}$. We get 
\begin{eqnarray} 
\begin{array}{ll} \displaystyle
\VF_{\geq 4, \geq 4}(x)=\frac{1-z^3}{1+z^3} \textrm{\LARGE{$($}} & \VF(x)-\VF(0)-x[x]\VF(x)+2x^2\VF(x)\left(\VF(x)-\VF(0)\right) \\ 
& \displaystyle +x^4\VF(x)^3-xz^5(\VF(x)-\VF(0)+x^2\VF(x)^2 -xz^2\VF(x))\\
& \displaystyle -2xz^2(\VF(x)-\VF(0)+x^2\VF(x)^2)-2\frac{x^2z^3}{1-z^3}\VF(x)(\VF(0)-z) ~\textrm{\LARGE{$)$}}.
\end{array} 
\label{eq:Vsup4sup4} \end{eqnarray}

Now, using Equations \Ref{eq:Rsup2toR} \Ref{eq:R2toR}  \Ref{eq:V3} \Ref{eq:decomposeV4} \Ref{eq:D4} and \Ref{eq:Vsup4sup4} we can replace  
$\VF_{\geq 2},~\VF_{2},~\VF_{3}$ and $\VF_{4}$ by their expression in Equation \Ref{eq:redecomposeV}. This establishes the following proposition.

\begin{prop} \label{prop:V} 
The generating function $\VF(x)=\VF(x,z)$ of non-separable near-triangulations in which all internal vertices have degree at least 5 satisfies:
\begin{eqnarray}
\hspace{-.2cm} \begin{array}{ll} 
\VF(x)=& \hspace{-.1cm}\displaystyle z+ 
\frac{1}{1+z^3}\pare{xz\VF(x)^2 +z\frac{\VF(x)-\VF_0}{x}-z^3\VF(x)} \vspace{.2cm} \\
& \hspace{-.3cm}\displaystyle -\frac{z^5(1-z^3)}{1+z^3}  \textrm{\Huge{$($}}  
\frac{\VF(x)-\VF_0-x\VF_1}{x^2} -z^2(2+z^3)\frac{\VF(x)-\VF_0}{x}-2\VF(x)(\VF_0-z) \vspace{.1cm}\\ 
& \displaystyle \hspace{1.6cm} +x^2\VF(x)^3-xz^2(2+z^3)\VF(x)^2 +2\VF(x)\left(\VF(x)-\VF_0\right) +z^7\VF(x) \textrm{\Huge{$)$}}
\end{array} \label{eq:V} 
\end{eqnarray} 
where $\VF_0=\VF(0)$ and  $\VF_1=[x]\VF(x)$ is the coefficient of $x$ in $\VF(x)$.

\end{prop}


\section{Algebraic equations for triangulations with high degree} \label{section:solution}
In the previous section, we have exhibited functional equations
concerning the families of \emph{near-triangulations} $\TT,\UU,\VV$. 
By definition, the generating functions $\TF(t),\UF(t),\VF(t)$ are power series in the main variable $z$ with polynomial coefficients in the secondary variable $x$.
We now solve these equations and establish algebraic equations for the families of \emph{triangulations} $\Fs,\GG,\HH$  in which vertices not incident to the root have degree at least $3,4,5$ respectively. As observed in Section \ref{section:preliminaires}, the generating functions $\FF(t),\GF(t),\HF(t)$ are closely related to the series $\TF(0),\UF(0),\VF(0)$ (see Equation \Ref{eq:alpha}).\\

Let us look at Equations \Ref{eq:T}, \Ref{eq:U} and \Ref{eq:V} satisfied by the series $\TF(x),~\UF(x)$ and $\VF(x)$ respectively.  We begin with Equation \Ref{eq:T}. This equation is (after multiplication by $x$) a polynomial equation in the main unknown series $\TF(x)$, the secondary unknown $\TF(0)$ and the variables $x,~z$. It is easily seen that this equation allows us to compute the coefficients of $\TF(x)$ (hence those of $\TF(0)$) iteratively. Moreover, we see by induction that the coefficients of this power series are polynomials in the secondary variable $x$. The same property holds for Equation \Ref{eq:U} (resp. \Ref{eq:V}): it defines the series $\UF(0)$ (resp. $\VF(0)$)  uniquely as a power series in $z$ with polynomial coefficients in $x$.\\

In some sense, Equations \Ref{eq:T}, \Ref{eq:U} and \Ref{eq:V} answer our enumeration problems. However, we want to \emph{solve} these equations, that is, to derive from them some equations for the series
$\TF(0),~\UF(0)$ and $\VF(0)$. Certain techniques for performing such manipulations appear in the combinatorics literature. 
In the cases of Equation  \Ref{eq:T} and \Ref{eq:U} which are quadratic in the main unknown series $\TF(x)$ and $\UF(x)$ we can routinely apply the so-called \emph{quadratic  method} \cite[Section 2.9]{Goulden:Lagrange}.  
This method allows one to solve polynomial equations which are quadratic in the bivariate unknown series and have one unknown univariate series. This method also applies to  Equation \Ref{eq:S(0)} concerning $\SF(x)$ and allows to prove Equation \Ref{eq:F}.
However, Equation \Ref{eq:V} concerning $\VF(x)$  is cubic in this series and involves two unknown univariate series ($\VF(0)$ and $[x]\VF(x)$). Very recently, Bousquet-Mélou and Jehanne proposed  a general method to solve polynomial equations of any degree in the bivariate unknown series and involving any number of unknown univariate series  \cite{MBM:quadratic}. We present their formalism.\\

Let us begin with Equation  \Ref{eq:T} concerning $\TF(0)$. We define the polynomial 
$$P(T,T_0,X,Z)=XZ+X^2ZT^2+ZT-ZT_0-XZ^3T-XT ~.$$
Equation  \Ref{eq:T} can be written as 
\begin{eqnarray} P(\TF(x),\TF(0),x,z)=0~.\label{eq:polyT} \end{eqnarray}
Let us consider the equation  $P_1'(\TF(x),\TF(0),x,z)=0$ , where $P_1'$ denotes the derivative of $P$ with respect to its first variable. This equation can be written as 
$$2x^2z\TF(x)+z-xz^3-x~=~0~.$$
This equation is not satisfied for a generic $x$. However, considered as an equation in $x$, it is straightforward to see that it admits a  unique power series solution $X(z)$.\\ 
Taking the derivative of Equation \Ref{eq:polyT} with respect to $x$ one obtains
$$\frac{\partial \TF(x)}{\partial x}\cdot P_1'(\TF(x),\TF(0),x,z)+ P_3'(\TF(x),\TF(0),x,z)=0,$$
where $P_3'$ denotes the derivative of $P$ with respect to its third variable. Substituting the series $X(z)$ for $x$ in that equation, we see that the series $X(z)$ is also a solution of the equation $P_3'(\TF(x),\TF(0),x,z)=0$. Hence, we have a system of three equations
$$
\begin{array}{rcl} 
P(\TF(X(z)),\TF(0),X(z),z)&=&0~, \\
P_1'(\TF(X(z)),\TF(0),X(z),z)&=&0~, \\
P_3'(\TF(X(z)),\TF(0),X(z),z)&=&0 ~,\\
\end{array} 
$$
for the three unknown series $\TF(X(z)),~\TF(0)$ and $X(z)$. This polynomial system can be solved by elimination techniques using either resultant calculations or Gr\"obner bases. Performing these eliminations one obtains an algebraic equation for $\TF(0)$:
$$\TF(0)=z-24{z}^{4}+3{z}^{7}+{z}^{10}+(32{z}^{3}+30{z}^{6}-4{z}^{9}-{z}^{12}) \TF(0) - 8z^2(1+z^3)^2 \TF(0)^2-16z^4 \TF(0)^3.$$
Using the fact that $\TF(0)=z+z\GF(z^3)$ we get the following theorem.

\begin{thm} \label{thm:G}
Let $\GG$ be the set of non-separable triangulations in which any
 vertex not incident to the root has degree at least 3, and let $\GF(t)$ be
 its generating function.
The series $\GF(t)$ is uniquely defined as a power series in $t$ by the algebraic equation:
\begin{eqnarray}
\begin{array}{l}
16t^{2}\GF(t)^{3}+8t(t^{2}+8t+1)\GF(t)^{2}\\
+(t^{4}+20t^{3}+50t^{2}-16t+1)\GF(t) +t^{2}(t^2+11t-1)=0. \label{eq:GFt} 
\end{array}
\end{eqnarray}

\end{thm}

Similar manipulations lead to a cubic equation for the set $\HH$.

\begin{thm}\label{thm:H}
Let $\HH$ be the set of  non-separable triangulations in which any vertex not incident to the root has degree at least 4, and let $\HF(t)$ be its generating function. 
The series $\HF(t)$ is uniquely defined as a power series in $t$ by the algebraic equation:
\begin{eqnarray}
\begin{array}{l}
16t^{2}(t-1)^{4}\HF(t)^{3}+(t^{8}+12t^{7}-14t^{6}-84t^{5}+207t^{4}-192t^{3}+86t^{2}-16t+1)\HF(t) \\
+8t(t-1)^{2}(t^{4}+4t^{3}-13t^{2}+8t+1)\HF(t)^{2}+t^{4}(t-1)(t^{3}+5t^{2}-8t+1)=0.
\end{array}
\label{eq:HFt} \end{eqnarray}
\end{thm}

For Equation \Ref{eq:V} concerning $\VF(0)$ the method is almost identical. We see that there is a polynomial $Q(V,V_0,V_1,x,z)$ such that Equation \Ref{eq:V} can be written as $Q(\VF(x),\VF(0),[x]\VF(x),x,z)=0$. But we can show that there are exactly \emph{two} series $X_1(z),~X_2(z)$ such that $Q_1'(\VF(X(z)),\VF(0),[x]\VF(x),X(z),z)=0$. Thus, we  obtain a system of 6 equations 
$$
\begin{array}{rcl} 
Q(\VF(X_i(z)),\VF(0),[x]\VF(x),X_i(z),z)&=&0 \\
Q_1'(\VF(X_i(z)),\VF(0),[x]\VF(x),X_i(z),z)&=&0\\
Q_3'(\VF(X_i(z)),\VF(0),[x]\VF(x),X_i(z),z)&=&0
\end{array}\hspace{1cm} i=1,~2 
$$
for the 6 unknown series   $\VF(X_1(z)),~\VF(X_2(z)),~X_1(z),~X_2(z),~\VF(0)$ and $[x]\VF(x)$.
This system can be solved via elimination techniques though the calculations involved are heavy. We obtain the following theorem.

\begin{thm}\label{thm:K}
Let $\KK$ be the set of  non-separable triangulations in which any vertex not incident to the root has degree at least 5, and let $\KF(t)$ be its  generating function. 
The series $\KF(t)$ is uniquely defined as a power series in $t$ by the algebraic equation:\\
\begin{eqnarray}
\sum_{i=0}^6 P_i(t) \KF(t)^i,
\label{eq:KFt} \end{eqnarray}
where the polynomials $P_i(t),i=0\ldots 6$ are given in Appendix \ref{appendixA}.
\end{thm}

\section{Constraining the vertices incident to the root}\label{section:constraining}
So far, we have established algebraic equations for the generating functions $\GF(t),\HF(t),\KF(t)$ of triangulations in which \emph{any vertex not incident to the root} has degree at least 3, 4, 5. The following theorems provide equations concerning  the generating functions $\GF^*(t),\HF^*(t)$ of triangulations in which \emph{any vertex} has degree at least 3, 4.

\begin{thm} \label{thm:G*}
Let $\GG^*$ be the set of  non-separable triangulations in which \emph{any} vertex has degree at least 3 and let $\GF^*(t)$ be its generating function. 
The series $\GF^*$ is related to the series $\GF$ of Theorem \ref{thm:G} by  
\begin{eqnarray} \GF^*(t)=(1-2t)\GF(t)~.\label{eq:G*G} \end{eqnarray}
\end{thm}

\begin{thm} \label{thm:H*}
Let $\HH^*$ be the set of  non-separable triangulations in which \emph{any} vertex  has degree at least 4 and let $\HF^*(t)$ be its generating function.   
The series $\HF^*$ is related to the series $\HF$ of Theorem \ref{thm:H} by  
\begin{eqnarray} \HF^*(t)=\frac{1-5t+5t^2-3t^3}{1-t}\HF(t)~.\label{eq:H*H} \end{eqnarray}
\end{thm}

\begin{figure}[hb!]
\begin{center}
\input{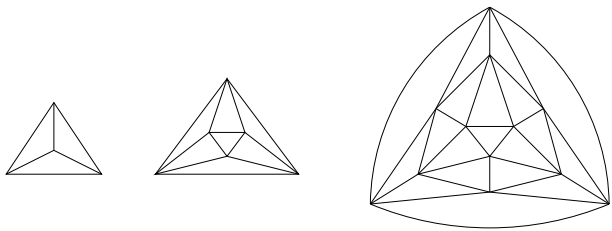}
\caption{The platonic solids: tetrahedron, octahedron, icosahedron.} \label{fig:platon}
\end{center}
\end{figure}

Let us make a few comments before proving these two theorems. First, observe that we can  deduce from Theorems \ref{thm:G} and \ref{thm:G*} (resp. \ref{thm:H} and \ref{thm:H*}) an algebraic equation for the generating function $\GF^*$ (resp. $\HF^*$) of triangulations in which \emph{any} vertex has degree at least 3 (resp. 4). The algebraic equation obtained for $\GF^*$ coincides with the result of Gao and Wormald \cite[Theorem 2]{Gao:triang_degresup}. From the algebraic equations we can routinely compute the first coefficients of our series: 
$$\hspace{-.2cm}\begin{array}{ll} & \GF^*(t)=t^2+3t^3+19t^4+128t^5+909t^6+6737t^7 +51683t^8+407802t^9+o(t^{9}), \\
 & \HF^*(t)= t^4+3t^5+12t^6+59t^7+325t^8+1875t^9+11029t^{10}+65607t^{11}+o(t^{11}).\\
\end{array}$$
Recall that the coefficient of $t^n$ in the series $\GF^*(t),~\HF^*(t)$ is the number of triangulations with $3n$ edges ($2n$ triangles, $n+2$ vertices) satisfying the required degree constraint. In the expansion of $\GF^*(t)$,  the smallest non-zero coefficient $t^2$ corresponds to the tetrahedron. In the expansion of $\HF^*(t)$, the smallest non-zero coefficient $t^4$ corresponds to the octahedron (see Figure \ref{fig:platon}). \\

We were unable to find an equation that would permit to count non-separable triangulations in which \emph{any vertex} has degree at least 5. However, we can use the algebraic equation \Ref{eq:KFt} to compute the first coefficients of the series $\KF(t)$:
$$\KF(t)={t}^{10}+8{t}^{11}+45{t}^{12}+209{t}^{13}+890{t}^{14}+3600{t}^{15}+14115{t}^{16}+54306{t}^{17}+o(t^{18}).$$
The first non-zero coefficient $t^{10}$ corresponds to the  icosahedron (see Figure \ref{fig:platon}).\\

\vspace{-.2cm} 
 \vspace{-.3cm}

In order to prove Theorems \ref{thm:G*} and \ref{thm:H*} we need some new notations. The set $\GG_{i,j,k}$ (resp. $\HH_{i,j,k}$) is the set of triangulations such that the root's origin  has degree $i$, the root's end has degree $j$, the third vertex of the root-face has degree $k$ and all internal vertices have degree at least 3 (resp. 4). For $\LL=\GG,\HH$ we define $\LL_{\geq i,j,k}= \bigcup_{l \geq i} \LL_{l,j,k}$ and with similar notation, $\LL_{\geq i,\geq j,k}$ etc. If $\LL$ is any of these sets, $\LF(t)$ is the corresponding generating function, where a map with $3n$ edges has contribution $t^n$.\\

\noindent \textbf{Proof of Theorem \ref{thm:G*}:} By definition, $\GG=\GG_{\geq 2, \geq 2, \geq 3}$ and $\GG^*=\GG_{\geq 3, \geq 3, \geq 3}$. Hence,
\begin{eqnarray} \GF^*(t)=\GF(t)-\GF_{2,2,\geq 3}(t)-2\GF_{2,\geq 3,\geq 3}(t).\label{eq:decomposeG} \end{eqnarray}
\ite The set $\GG_{2,2,\geq 3}$ is empty, hence $\GF_{2,2,\geq 3}(t)=0$. \\
\ite The set $\GG_{2,\geq 3,\geq 3}$ is in bijection with $\GG_{\geq 1,\geq 1,\geq 3}=\GG$ by the mapping represented in Figure \ref{fig:G2sup3toG}. This gives $\GF_{2,\geq 3,\geq 3}(t)=t\GF(t)$. \\
Plugging these results in \Ref{eq:decomposeG} proves the theorem. \findem

\begin{figure}[ht!]
\begin{center}
\input{G2sup3toG.pstex_t}
\caption{A bijection between $\GG_{2, \geq 3,\geq 3}$ and $\GG$ (resp. $\HH_{2, \geq 3, \geq 3}$ and $\HH$).} \label{fig:G2sup3toG}
\end{center}
\end{figure}

\begin{figure}[ht!]
\begin{center}
\input{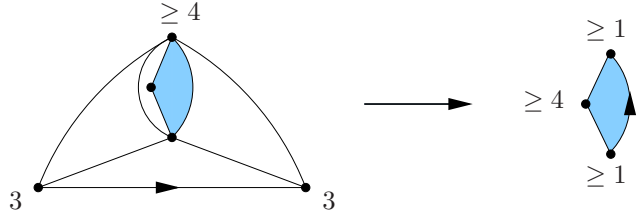}
\caption{A bijection between $\HH_{3,3,\geq 4}$ and $\HH$.} \label{fig:H33toH}
\end{center}
\end{figure}

\noindent \textbf{Proof of Theorem \ref{thm:H*}:}
By definition, $\HH^*=\HH_{\geq 4, \geq 4, \geq 4}$. Hence, 
\begin{eqnarray} \HF^*(t)=\HF_{\geq 3, \geq 3, \geq 4}(t)-\HF_{3,3,\geq 4}(t)-2\HF_{3,\geq 4,\geq 4}(t). \label{eq:decomposeH} \end{eqnarray}
Recall that  $\HH=\HH_{\geq 1, \geq 1, \geq 4}=\HH_{\geq 2, \geq 2, \geq 4}$.\\
\ite Clearly, $\HF_{\geq 3, \geq 3, \geq 4}(t)= \HH_{\geq 2,\geq 2,\geq 4}(t) -\HH_{2,2,\geq 4}(t)-2\HH_{2,\geq 3,\geq 3}(t)$. \\
\ite The set $\HF_{2,2,\geq 4}(t)$ is empty, hence $\HH_{2,2,\geq 4}(t)=0$. \\
\ite The set $\HH_{2,\geq 3,\geq 3}$ is in bijection with $\HH_{\geq 1,\geq 1,\geq 4}=\HH$ by the mapping represented in Figure \ref{fig:G2sup3toG}, hence $\HH_{2,\geq 3,\geq 3}(t)=t\HF(t)$. \\
This gives 
\begin{eqnarray}\HF_{\geq 3,\geq 3,\geq 4}(t)=(1-2t)\HF(t).\label{eq:Hsup3sup3sup4} \end{eqnarray} 
\ite The set $\HH_{3,3,\geq 4}$ is in bijection with $\HH_{\geq 1,\geq 1,\geq 4}=\HH$ by the mapping represented in Figure \ref{fig:H33toH}. This gives 
\begin{eqnarray}\HF_{3,3,\geq 4}(t)=t^2\HF(t).\label{eq:H33sup4}\end{eqnarray}

\begin{figure}[ht!]
\begin{center}
\input{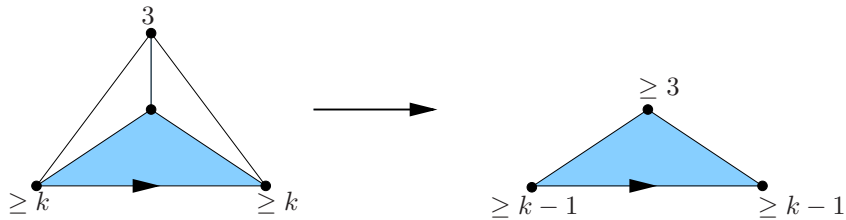}
\caption{A bijection between $\HH_{\geq k,\geq k,3}$ and   $\HH_{\geq k-1,\geq k-1,\geq 3}$.} \label{fig:Hsupksupk3}
\end{center}
\end{figure}

\ite For any integer $k$ greater than 2, the set $\HH_{\geq k,\geq k,3}$ is in bijection with the set   $\HH_{\geq k-1,\geq k-1,\geq 3}$ by the mapping represented in Figure \ref{fig:Hsupksupk3}. This gives
\begin{eqnarray}  \HF_{\geq k,\geq k,3}(t)=t\HF_{\geq k-1,\geq k-1,\geq 3}(t)~~~\textrm{ for all } k\geq2. \label{eq:relation-k} \end{eqnarray}

\noindent Using Equation \Ref{eq:relation-k} for $k=4$ and then for $k=3$ (and trivial symmetry properties), we get \\
$\hspace{-.2cm} 
\begin{array}{ll}
\HF_{3,\geq 4,\geq 4}(t) & =~\HF_{\geq 4,\geq 4,3}(t)~=~t\HF_{\geq 3,\geq 3,\geq 3}(t)~=~t\HF_{\geq 3,\geq 3,\geq 4}(t)+t\HF_{\geq 3,\geq 3,3}(t) \\
& =~t\HF_{\geq 3,\geq 3,\geq 4}(t)+t^2\HF_{\geq 2,\geq 2,\geq 3}(t).
\end{array}$\\
\iten By Equation  \Ref{eq:Hsup3sup3sup4}, we have  $\HF_{\geq 3,\geq 3,\geq 4}(t)=(1-2t)\HF(t)$. \\
\iten Using Equation \Ref{eq:relation-k} for $k=2$  gives  
$$\HF_{\geq 2,\geq 2,\geq 3}(t)=\HF_{\geq 2,\geq 2,\geq 4}(t)+\HF_{\geq 2,\geq 2,3}(t)=\HF(t)+t\HF_{\geq 1,\geq 1,\geq 3}(t).$$ 
Given that $\HH_{\geq 1,\geq 1,\geq 3}=\HH_{\geq 2,\geq 2,\geq 3}$, we get $\displaystyle \HF_{\geq 2,\geq 2,\geq 3}(t)=\frac{1}{1-t}  \HF(t)$.\\ 
Thus, we obtain  
\begin{eqnarray} \HF_{3,\geq 4,\geq 4}(t)=\frac{t(1-2t+2t^2)}{1-t}\HF(t). \label{eq:H3sup4sup4}\end{eqnarray}

\noindent Plugging Equations \Ref{eq:Hsup3sup3sup4}, \Ref{eq:H33sup4} and \Ref{eq:H3sup4sup4} in Equation \Ref{eq:decomposeH} proves the theorem.
\findem


\section{Asymptotics} \label{section:asymptotic}
In Section \ref{section:solution}, we established algebraic equations for the generating functions $\LF=\FF,\GF,\HF,\KF$ of  non-separable triangulations in which any vertex not incident to the root has degree at least $d=2,3,4,5$ (Equations  \Ref{eq:F}, \Ref{eq:GFt}, \Ref{eq:HFt} and \Ref{eq:KFt}). We will now derive the asymptotic form of the number  $l_n=f_n,g_n,h_n,k_n$  of maps with  $3n$ edges in each family by analyzing the singularities of the generating function $\LF=\FF,\GF,\HF,\KF$ ($l_n$ is the coefficient of $t^n$ in $\LF$).  The principle of this method is a general correspondence between the expansion of a generating function at its dominant singularities and the asymptotic form of its coefficients \cite{Flajol:Sing1,Flajolet:analytic}.\\
\begin{lemma} \label{thm:sugularexpansion}
Each of the generating functions $\LF=\FF,\GF,\HF,\KF$ has a unique dominant singularity $\rho_L>0$ and a singular expansion with singular exponent $\frac{3}{2}$ at $\rho_L$, in the sense that 
\begin{eqnarray}
\LF(t)=\alpha_L+\beta_L(1-\frac{t}{\rho_L})+\gamma_L (1-\frac{t}{\rho_L})^{3/2}+O((1-\frac{t}{\rho_L})^{2}),
\end{eqnarray}
with $\gamma_L\neq 0$. The dominant singularities of the series $\FF$ and $\GF$ are respectively $\rho_F=\frac{2}{27}$ and $\rho_G=\frac{3\sqrt{3}-5}{2} $. The dominant singularities $\rho_H$ and $\rho_K$ of the series $\HF$ and $\KF$ are defined by algebraic equations given in Appendix \ref{appendixB}. 
\end{lemma}

\textbf{Proof (sketch):}
The (systematic) method we follow is described in  \cite[Chapter VII.4]{Flajolet:analytic}). Calculations were performed using the Maple package \emph{gfun} \cite{Salvy:gfun}.\\
Let us denote  generically by $\rho_L$  the radius of convergence of the series $\LF$ and by  $Q(\LF,t)$ the algebraic equation satisfied by $\LF$ (Equations \Ref{eq:F}, \Ref{eq:GFt}, \Ref{eq:HFt} and \Ref{eq:KFt}).
It is known that the singular points of the series $\LF$ are among the roots of the polynomial $R(t)=D(t)\Delta(t)$ where $D(t)$ is the dominant coefficient of $Q(y,t)$ and  $\Delta(t)$ is the  discriminant of $Q(y,t)$  considered as a polynomial in $y$.  Moreover, since the series $\LF$ has non-negative coefficients, we know (by Pringsheim's Theorem) that the point $t=\rho_L$ is singular. In our cases, the smallest positive root of $R(t)$ is found to be indeed a singular point of the series $\LF$. (This requires to solve some connection problems that we do not detail.)  Moreover, no other root of $R(t)$ has the same modulus. This proves that the series $\LF$ has a unique dominant singularity.\\
The second step is to expand the series $\LF$ near its singularity $\rho_L$. This calculation can be performed using \emph{Newton's polygon method}  (see \cite[Chapter VII.4]{Flajolet:analytic}) which is implemented in the \emph{algeqtoseries} Maple command \cite{Salvy:gfun}. \findem

From Lemma \ref{thm:sugularexpansion}, we can deduce the  asymptotic form of the  number $l_n=f_n,g_n,h_n,k_n$ of  non-separable triangulations of size $n$ in each family.
\begin{thm}
The number $l_n=f_n,g_n,h_n,k_n$ of  non-separable triangulations of size $n$ ($3n$ edges) in which any vertex not incident to the root has degree at least $d=2,3,4,5$ has asymptotic form 
$$l_n~ \sim \lambda_L n^{-5/2} \parfrac{1}{\rho_L}^n . $$
The growth constants $\rho_F,\rho_G,\rho_H,\rho_K$ are given in Lemma \ref{thm:sugularexpansion}. Numerically, 
$$\frac{1}{ \rho_F}=13.5,~~\frac{1}{ \rho_G}\approx 10.20,~~\frac{1}{ \rho_H}\approx 7.03,~~\frac{1}{ \rho_K}\approx 4.06~.$$
\end{thm}

\noindent \textbf{Remark:} The subexponential factor $n^{-5/2}$ is typical of planar maps families (see for instance  \cite{Banderier:Airy} where 15 classical families of maps are listed all displaying this subexponential factor $n^{-5/2}$). \\

\noindent \textbf{Remark:} Using Theorems \ref{thm:G*} and \ref{thm:H*}, it is easily seen that the series $\LF^*=\GF^*,\HF^*$ has dominant singularity  $\rho_L=\rho_G,\rho_H$ with singular exponent $\frac{3}{2}$ at $\rho_L$:
$$ \LF(t)=\alpha_L^*+\beta_L^*(1-\frac{t}{\rho_L})+\gamma_L^* (1-\frac{t}{\rho_L})^{3/2}+O((1-\frac{t}{\rho_L})^{2}).$$
Therefore, we obtain the asymptotic form 
$$l_n^*~ \sim \lambda_L^* n^{-5/2} \parfrac{1}{\rho_L}^n  $$
for the number $l_n^*=g_n^*,h_n^*$ of  non-separable triangulations of size $n$ with vertex degree at least $d=3,4$. Hence, the numbers $l_n^*$ and $l_n$ are equivalent up to a (known) constant multiplicative factor $\displaystyle \frac{\lambda_L^*}{\lambda_L}$:
$$
\begin{array}{l}
\displaystyle \frac{\lambda_G^*}{\lambda_G}~=~\frac{\gamma_G^*}{\gamma_G}~=~1-2\rho_G=6-3\sqrt{3},\vspace{.1cm}\\
\displaystyle \frac{\lambda_H^*}{\lambda_H}~=~\frac{\gamma_H^*}{\gamma_H}~=~\frac{1-5\rho_H+5{\rho_H}^2-3{\rho_H}^3}{1-\rho_H}.
\end{array}
$$
We do not have such precise information about the asymptotic form of the number $k_n^*$ of  non-separable triangulations of size $n$ ($3n$ edges) with vertex degree at least 5. However, we do know that  $\displaystyle k_n^*=\Theta(k_n)=\Theta(n^{-5/2} {\rho_K}^{-n}  )$. Indeed, we clearly have $k_n^*\leq k_n$ and, in addition, $k_n^*\geq k_{n-9}\sim {\rho_K}^9 k_n$. The latter inequality is  proved by observing that the operation of replacing the root-face of a triangulation by an icosahedron is an injection from the set of triangulations of size $n$ in which any vertex \emph{not incident to the root} has degree $5$ to the set of triangulations of size $n+9$ in which \emph{any} vertex has degree at least 5.\\

\section{Concluding remarks}
We have established algebraic equations for the generating functions of non-separable triangulations in which \emph{any vertex not incident to the root} has degree at least $d=3, 4, 5$. We have also established algebraic equations for non-separable triangulations in which \emph{any vertex} has degree at least $d=3, 4$. However, have not found a similar result for $d=5$. The algebraic equations we have obtained can be converted into differential equations (using for instance the \emph{algeqtodiffeq} Maple command available in the \emph{gfun} package \cite{Salvy:gfun}) from which one can compute the coefficients of the series in a linear number of operations.  Moreover, the asymptotic form of their coefficients can also be found routinely from the algebraic equations.\\    

The approach we have adopted is based on a classic decomposition scheme allied with a generating function approach. Alternatively, it is possible to obtain some of our results by a compositional approach. This is precisely the method followed by Gao and Wormald to obtain the algebraic equation concerning non-separable triangulations in which any vertex has degree at least 3  \cite{Gao:triang_degresup}. We were also able to obtain  the algebraic equation concerning non-separable triangulations in which any vertex has degree at least 4 by this composition method \cite{OB:pipeau-triangulations-haut-degres-composition}. However, we do not see how to apply this method to non-separable triangulations in which vertices not incident to the root have degree at least 5.\\

Recently, Poulalhon and Schaeffer gave a bijective proof  based on the \emph{conjugacy classes of tree} for the number of non-separable triangulations \cite{Schaeffer:triangulation}. However, it is dubious that this approach should apply for the families $\HF,~\KF$ of non-separable triangulations in which vertices have degree at least $d=4,5$. Indeed, for a large number of families of maps $\LL$, the generating function $\LF(t)$ is \emph{Lagrangean}, that is, there exists a series $\XF(t)$ and two rational functions $\Psi,\Phi$ satisfying 
 $$\LF(t)=\Psi(\XF(t)) ~~\textrm{and}~~ \XF(t)=t\Phi(\XF(t))$$ 
(see for instance \cite{Banderier:Airy} where 15 classical families are listed together with a Lagrangean parametrization). Often, a parametrization can be found such that the series $\XF(t)$ looks like the generating function of a family of trees (i.e. $\Phi(x)$ is a series with non-negative coefficients) suggesting that a bijective approach exists based on the enumeration of certain trees \cite{MBM:Ising2,DiFrancesco:census,DiFrancesco:hardbicubic}.
However, it is known that an algebraic series is Lagrangean if and only if the genus of the algebraic equation is 0 \cite[Chapter 15]{Abhyankar:algebraic}. In our case, the algebraic equations defining the series $\FF$, $\GF,~\HF$ and $\KF$ have respective genus 0, 0, 2 and 25. (The genus can be computed using the Maple command \emph{genus}.) Thus, whereas the series $\FF$, $\GF$ are Lagrangean (with a parametrization given in Appendix \ref{appendixC}), the series $\HF,~\KF$ are not. \\



Lastly, we claim some generality to our approach. Here, we have focused on non-separable triangulations, but it is possible to practice the same kind of manipulations for \emph{general} triangulations. The method should also apply to some other families of maps, like quadrangulations. Thus, a whole new class of map families is expected to have algebraic generating functions. \\

\noindent \textbf{Acknowledgments:} I would like to thank Mireille Bousquet-Mélou. This paper and I benefited greatly from her suggestions and support.

\bibliography{allref}
\bibliographystyle{plain}

\newpage
\appendix
\section{Coefficients of the algebraic equation \Ref{eq:KFt}} \label{appendixA}
The coefficients $P_i(t),i=0..6$ in the algebraic equation \Ref{eq:KFt} are:\\

\noindent $P_0(t)=t^{10}(-1+82552t^{11}-163081t^{12}+277796t^{13}-308156t^{14}-443851t^{16}+t^{34}+13t+32t^{31}+454t^{5}-2434t^{6}-5762t^{8}+4373t^{7}-53961t^{10}+23037t^{9}+354387t^{15}+163964t^{20}-28454t^{21}-38408t^{22}+36713t^{23}-11737t^{24}+t^{33}+2t^{32}-278t^{25}+242t^{28}-1678t^{27}+2714t^{26}+36t^{29}-64t^{30}-70t^{2}+180t^{3}-195t^{4}-273662t^{19}+122688t^{18}+262614t^{17})$,\\

\noindent $P_1(t)=(1-594873t^{11}+1078572t^{12}-1457943t^{13}+1921912t^{14}+1327736t^{16}+1462t^{38}-3168t^{37}-611t^{39}+25956t^{35}-56515t^{34}-3826t^{36}-21t-467567t^{31}-4545t^{5}+3916t^{6}+60304t^{8}-13364t^{7}+275068t^{10}-142715t^{9}-2t^{42}+9t^{43}+t^{44}-2338117t^{15}-4673450t^{20}+5167054t^{21}-1145738t^{22}-2425736{t}^{23}+2298353t^{24}+66635t^{33}+90827t^{32}+559893t^{25}-874518t^{28}+2995671t^{27}-3225500t^{26}-526335t^{29}+763474t^{30}+68t^{41}+75t^{40}+193t^{2}-988t^{3}+2913t^{4}+1719643t^{19}-945302t^{18}+541155t^{17})$,\\

\noindent $P_2(t)= t(8+2011979t^{11}-1422607t^{12}+2174211t^{13}-4910332t^{14}-9095603{t}^{16}-814t^{38}+688t^{37}+306t^{39}-16997t^{35}+43703t^{34}+1292t^{36}-4t+370239t^{31}-3000t^{5}+20421t^{6}-268574t^{8}+72382t^{7}-1309172t^{10}+527412t^{9}+8t^{42}+5383141t^{15}+31153077t^{20}-16211612{t}^{21}-2143067t^{22}+7886923t^{23}-2902691t^{24}-50536{t}^{33}-26161t^{32}-4609909t^{25}+156674t^{28}-3199107{t}^{27}+6488106t^{26}+970079t^{29}-902321t^{30}+12t^{41}+4t^{40}-556t^{2}+3851t^{3}-8840t^{4}-18494688t^{19}-9439987t^{18}+17752182t^{17})  $,\\

\noindent $P_3(t)= t^{2}(16+1278321t^{11}-2978655t^{12}+1697247t^{13}+5975715t^{14}+54631824t^{16}+166t^{38}-90t^{37}-32t^{39}+3984t^{35}-13104t^{34}-868t^{36}-192t-105251t^{31}+17247t^{5}-36981t^{6}+521925t^{8}-74982t^{7}+835782t^{10}-1142394t^{9}-29427957t^{15}-39935486t^{20}+7773505t^{21}+6824437t^{22}-5541795t^{23}-1619262t^{24}+18648t^{33}+4941t^{32}+5146785t^{25}+349680t^{28}+880004t^{27}-3411645t^{26}-600239t^{29}+358687t^{30}+16t^{40}+1046t^{2}-2554t^{3}-397t^{4}+60017232t^{19}-26467945t^{18}-34977363t^{17})  $,\\

\noindent $P_4(t)= 9t^{5}(t-1)^{2} (8+722739t^{11}-1888278t^{12}+1483343t^{13}+679876t^{14}+1099122t^{16}-84t-20t^{31}+9250t^{5}-17908t^{6}+144652t^{8}-22565t^{7}+87721t^{10}-234335t^{9}-1820089t^{15}-5409t^{20}-64607t^{21}+41918t^{22}-12628t^{23}-1362t^{24}+8t^{32}+6200t^{25}-189t^{28}+1127t^{27}-3809t^{26}-103t^{29}+84t^{30}+368t^{2}-583t^{3}-2069t^{4}+110521t^{19}-69119t^{18}-243772t^{17})  $,\\

\noindent $P_5(t)= 81t^{8}(t-1)^{4}(1+25926t^{11}-14080t^{12}+2973t^{13}-369t^{14}+348t^{16}-9t+2118t^{5}-2936t^{6}+23913t^{8}-4134t^{7}-6330t^{10}-25946t^{9}-970t^{15}+12t^{20}-22t^{21}+12t^{22}-3t^{23}+t^{24}+30t^{2}-15t^{3}-747t^{4}+42t^{19}-219t^{18}+405t^{17}
)   $,\\

\noindent $P_6(t)=59049t^{15}(t+1) (t-1)^{9}$.\\

\section{Algebraic equations for the dominant singularity of the series $\HF(t)$ and $\KF(t)$}\label{appendixB}
The dominant singularity $\rho_H$ (resp. $\rho_K$) of the generating function $\HF(t)$ (resp. $\KF(t)$) is the smallest positive root of the polynomial $r_H(t)$ (resp. $r_K(t)$) where \\
$r_H(t)=2-17t+22t^2-10t^3+2t^4, $\\
 and\\
$r_K(t)=256-5504t+51744{t}^{2}-265664{t}^{3}+755040{t}^{4}-1069751{t}^{5}+1411392{t}^{6}-9094370{t}^{7}+30208920{t}^{8}-14854607{t}^{9}-106655904{t}^{10}+169679596{t}^{11}+1693392{t}^{12}+58535932{t}^{13}-263701752{t}^{14}-751005332{t}^{15}+2215033200{t}^{16}-2276240390{t}^{17}+2301677920{t}^{18}-1558097344{t}^{19}-2448410184{t}^{20}+6223947236{t}^{21}-7440131352{t}^{22}+6100648148{t}^{23}+1602052848{t}^{24}-9604816702{t}^{25}+6144202392{t}^{26}+996698032{t}^{27}+551560496{t}^{28}-3299013583{t}^{29}-728097928{t}^{30}+4881643814{t}^{31}-3845803168{t}^{32}+494467523{t}^{33}+1677669800{t}^{34}-1787552140{t}^{35}+825330824{t}^{36}+1529759{t}^{37}-340280968{t}^{38}+301075034{t}^{39}-121555768{t}^{40}-1710967{t}^{41}+37850432{t}^{42}-27659392{t}^{43}+9430688{t}^{44}-152352{t}^{45}-1901664{t}^{46}+1245152{t}^{47}-400416{t}^{48}+47744{t}^{49}+30720{t}^{50}-22528{t}^{51}+7680{t}^{52}-1792{t}^{53}+256{t}^{54}.$\\

\section{Lagrangean parametrization for the series $\FF(t)$, $\GF(t)$ and $\GF^*(t)$}\label{appendixC}
The series $\FF(t)$ has the following Lagrangean parametrization:
$$\FF(t)=\frac{\XF(1+\XF)}{2},$$
where
$$\XF\equiv \XF(t)=2t(1+\XF(t))^3.$$

The series $\GF(t)$ and $\GF^*(t)$ have the following Lagrangean parametrization:
$$
\begin{array}{lll}
\GF(t)&=&\displaystyle 2t\YF(1+ \YF) (1-\YF-{\YF}^{2}),\vspace{.2cm}\\
\GF^*(t)&=&\displaystyle 4t^2 (1+ \YF) (1-\YF-{\YF}^{2})(1+3\YF+6{\YF}^{2}+2{\YF}^{3}),
\end{array}
$$
where
$$\YF\equiv \YF(t)=2t(1+ \YF(t))  (1+4\YF(t)+2\YF(t)^{2}).$$ 

\end{document}